\documentclass[12pt,twoside,reqno,leqno]{amsart}
\usepackage{url}
\usepackage{amsmath}
\usepackage{amssymb}
\usepackage{amsfonts}
\usepackage{amsthm}
\usepackage[mathscr]{eucal}
\usepackage{amscd}
\usepackage{amsmath}
\usepackage{latexsym}
\usepackage[dvips]{graphics}
\usepackage{graphics,epstopdf}
\usepackage[pdftex]{graphicx}
\theoremstyle{plain}

\theoremstyle{definition}

 \theoremstyle{remark}
 \newtheorem{rem}{Remark}
\newtheorem{example}{Example}

\allowdisplaybreaks

\begin{document}
\title[Relation-Theoretic Metrical Coincidence Theorems] {Relation-Theoretic Metrical Coincidence Theorems}
\author[Aftab Alam and Mohammad Imdad]{Aftab Alam$^{\ast}$ and Mohammad Imdad}
 \thanks{$^\ast$Correspondence:\;{\rm aafu.amu@gmail.com}\bigskip}
\maketitle
\begin{center}
{\footnotesize Department of Mathematics\\Aligarh Muslim
University\\ Aligarh-202002, India.\\
Email addresses: aafu.amu@gmail.com, mhimdad@gmail.com\\}
\end{center}
{\footnotesize{\noindent {\bf Abstract.}}  In this article, we
generalize some frequently used metrical notions such as:
completeness, closedness, continuity, $g$-continuity and compatibility to
relation-theoretic setting and utilize these relatively
weaker notions to prove results on the existence and uniqueness of coincidence
points involving a pair of mappings defined on a metric space endowed with an
arbitrary binary relation. Particularly, under universal relation
our results deduce the classical coincidence point theorems of
Goebel (Bull. Acad. Pol. Sci. S$\acute{\rm e}$r. Sci. Math. Astron.
Phys. 16 (1968) 733-735) and Jungck (Int. J. Math. Math. Sci. 9 (4)
(1986) 771-779). In process our results generalize, extend, modify and unify several well-known results
especially those obtained in Alam and Imdad (J. Fixed Point Theory Appl. 17 (4) (2015) 693-702), Karapinar $et\;al.$ (Fixed
Point Theory Appl. 2014:92 (2014) 16 pp), Alam $et\;al.$ (Fixed
Point Theory Appl. 2014:216 (2014) 30 pp), Alam and Imdad (Fixed
Point Theory, in press) and Berzig (J. Fixed Point Theory Appl. 12
(1-2) (2012) 221-238.\\


}

\section{Introduction and Preliminaries}
\label{SC:Introduction and Preliminaries}

Throughout this manuscript, $\mathbb{N}$, $\mathbb{N}_0$, $\mathbb{Q}$ and $\mathbb{R}$
denote the sets of natural numbers, whole
numbers, rational
numbers and real numbers respectively ($i.e.$
$\mathbb{N}_0=\mathbb{N} \cup \{0
\}$). For the sake
of completeness, firstly we recall some known relevant definitions.\\

\noindent\textbf{Definition 1 \cite{CP1,CP2}.} Let $X$ be a nonempty
set and $f$ and $g$ two self-mappings on $X$. Then
\begin{enumerate}
\item[{(i)}] an element $x\in X$ is called a
coincidence point of $f$ and $g$ if
$$g(x)=f(x),$$
\item[{(ii)}] if $x\in X$ is a coincidence point of $f$ and $g$
and $\overline{x}\in X$ such that $\overline{x}=g(x)=f(x)$, then
$\overline{x}$ is called a point of coincidence of $f$ and $g$,
\item[{(iii)}] if $x\in X$ is a coincidence point of $f$ and $g$
such that $x=g(x)=f(x)$, then $x$ is called a common fixed point of
$f$ and $g$,
\item[{(iv)}] $f$ and $g$ are called commuting if
$$g(fx)=f(gx)\hspace{0.5cm}\forall~x\in X\;{\rm and}$$
\item[{(v)}] $f$ and $g$ are called weakly compatible (or partially commuting or
coincidentally commuting) if $f$ and $g$ commute at their
coincidence points, $i.e.$, for any $x\in X,$
$$g(x)=f(x)\Rightarrow g(fx)=f(gx).$$
\end{enumerate}

\noindent\textbf{Definition 2 \cite{CP3,CP4,g-C}.} Let $(X,d)$ be a
metric space and $f$ and $g$ two self-mappings on $X$. Then
\begin{enumerate}
\item[{(i)}] $f$ and $g$ are called weakly commuting if for all $x\in X,$
$$d(gfx,fgx)\leq d(gx,fx)\;{\rm and}$$
\item[{(ii)}] $f$ and $g$ are called compatible if
$$\lim\limits_{n\to \infty}d(gfx_n,fgx_n)=0$$
whenever $\{x_n\}$ is a sequence in $X$ such that
$$\lim\limits_{n\to \infty}g(x_n)=\lim\limits_{n\to \infty}f(x_n),$$
\item[{(iii)}] $f$ is called $g$-continuous at some $x\in X$ if for all sequences $\{x_n\}\subset X$,
$$g(x_n)\stackrel{d}{\longrightarrow} g(x)\Rightarrow f(x_n)\stackrel{d}{\longrightarrow} f(x).$$
Moreover, $f$ is called $g$-continuous if it is $g$-continuous at
each point of $X$.
\end{enumerate}

Recently, Alam and Imdad \cite{RT1} extended the classical Banach
contraction principle to complete metric space endowed with a binary
relation and observed that the partial order (see Nieto and
Rodr\'{\i}guez-L\'{o}pez \cite{PF2}), preorder (see Turinici
\cite{T-LCP}), transitive relation (see Ben-El-Mechaiekh
\cite{BR3}), tolerance (see Turinici \cite{T-RRF,T-NLF}), strict
order (see Ghods $et\;al.$ \cite{C-SO}), symmetric closure (see
Samet and Turinici \cite{BR1}) $etc$ utilized in several well-known
metrical fixed point theorems are not necessary and we can extend such
results for an arbitrary binary relation. In this context, the
contraction condition is relatively weaker than usual contraction as
it is required to hold for only those elements which are related in the
underlying
relation rather than for every pair of elements.\\

The aim of this paper is to extend our results proved in \cite{RT1}
to prove some existence and uniqueness results on coincidence points
in metric space endowed with an arbitrary binary relation. In proving the results, we employ our newly introduced notions such as:
$\mathcal{R}$-completeness, $\mathcal{R}$-closedness, $\mathcal{R}$-continuity, $(g,\mathcal{R})$-continuity,
$\mathcal{R}$-compatibility, $\mathcal{R}$-connected sets $etc.$. In this
course, we also observe that our results combine the idea contained
in Karapinar $et\;al.$ \cite{FIC9} as the set {\it M} (utilized by
Karapinar $et\;al.$ \cite{FIC9}) being subset of $X^2$ is infact a
binary relation on $X$. As consequences of our
newly proved results, we deduce several other established metrical coincidence point theorems. Finally, we furnish some
illustrative examples to demonstrate our results.

\section{Relation-Theoretic Notions and Auxiliary Results}
\label{SC:Relation-Theoretic Notions and Auxiliary  Results}

In this section, to make our exposition self-contained, we give some
definitions and basic results related to our main results.\\

\noindent\textbf{Definition 3 \cite{ST}.} Let $X$ be a nonempty set.
A subset $\mathcal{R}$ of $X^2$ is called a binary
relation on $X$.
Notice that for each pair $x,y\in X$,
one of the following holds:\\
\indent\hspace{0.5cm}(i) $(x,y)\in \mathcal{R}$; means
that ``$x$ is $\mathcal{R}$-related to $y$" or ``$x$ relates to
$y$ under
$\mathcal{R}$". Sometimes, we write $x\mathcal{R}y$ instead of $(x,y)\in\mathcal{R}$.\\
\indent\hspace{0.5cm}(ii) $(x,y)\not\in \mathcal{R}$; means
that ``$x$ is not $\mathcal{R}$-related to $y$" or ``$x$ doesn't
relate to $y$ under
$\mathcal{R}$".\\

Trivially, $X^2$ and $\emptyset$ being subsets of $X^2$ are binary
relations on $X$, which are respectively called the universal
relation (or full relation) and empty relation.\\

Throughout this paper, $\mathcal{R}$ stands for a nonempty binary
relation but for the sake of simplicity, we often write `binary
relation' instead of `nonempty binary relation'.\\

\noindent\textbf{Definition 4 \cite{RT1}.} Let $\mathcal{R}$ be a
binary relation on a nonempty set $X$ and $x,y\in X$. We say
that $x$ and $y$ are $\mathcal{R}$-comparative if either $(x,y)\in
\mathcal{R}$ or $(y,x)\in \mathcal{R}$. We denote it by $[x,y]\in
\mathcal{R}$.\\

\noindent{\bf Proposition 1.} If $(X,d)$ is a metric space,
$\mathcal{R}$ is a binary relation on $X$, $f$ and $g$ are two
self-mappings on $X$ and $\alpha\in [0,1)$, then the following contractivity conditions are equivalent:\\
\indent\hspace{0.1cm}(I) $d(fx,fy)\leq\alpha d(gx,gy)\;\;\forall~ x,y\in X$ with $(gx,gy)\in \mathcal{R}$,\\
\indent\hspace{0.1cm}(II) $d(fx,fy)\leq\alpha d(gx,gy)\;\;\forall~ x,y\in X$ with $[gx,gy]\in \mathcal{R}$.\\
{\noindent{Proof.}} The implication (II)$\Rightarrow$(I) is trivial.
On the other hand, suppose that (I) holds. Take $x,y\in X$ with
$[gx,gy]\in \mathcal{R}$. If $(gx,gy)\in \mathcal{R}$, then (II) is
directly follow from (I). Otherwise, in case $(gy,gx)\in
\mathcal{R}$, using symmetry of $d$ and (I), we obtain
$$d(fx,fy)=d(fy,fx)\leq\alpha d(gy,gx)=\alpha d(gx,gy).$$ This proves
that (I)$\Rightarrow$(II).\\

\noindent\textbf{Definition 5 \cite{ST,RA}.} A binary relation $\mathcal{R}$ on a nonempty set $X$ is called\\
\indent\hspace{0.5cm}{\em reflexive} if $(x,x)\in \mathcal{R}~~~\forall x\in X$,\\
\indent\hspace{0.5cm}{\em symmetric} if whenever $(x,y)\in \mathcal{R}$ then $(y,x)\in \mathcal{R}$,\\
\indent\hspace{0.5cm}{\em antisymmetric} if whenever $(x,y)\in \mathcal{R}$ and $(y,x)\in \mathcal{R}$ then $x=y$,\\
\indent\hspace{0.5cm}{\em transitive} if whenever $(x,y)\in \mathcal{R}$ and $(y,z)\in \mathcal{R}$ then $(x,z)\in \mathcal{R}$,\\
\indent\hspace{0.5cm}$\bullet$\;{\em complete} or {\em connected} or {\em dichotomous} if $[x,y]\in \mathcal{R}~~~\forall~ x,y\in X$,\\
\indent\hspace{0.5cm}$\bullet$\;{\em weakly complete} or {\em weakly
connected} or {\em trichotomous} if $[x,y]\in \mathcal{R}$\\
\indent\hspace{0.7cm} or $x=y~~~\forall~ x,y\in X$.\\

\noindent\textbf{Definition 6 \cite{BR1,ST,RA,BR01,BR02,BR03}.} A binary
relation
$\mathcal{R}$ defined on a nonempty set $X$ is called\\
\indent\hspace{0.5cm}$\bullet$\;{\em amorphous} if $\mathcal{R}$ has no specific properties at all,\\
\indent\hspace{0.5cm}$\bullet$\;{\em strict order} or {\em sharp
order} if
$\mathcal{R}$ is irreflexive and transitive,\\
\indent\hspace{0.5cm}$\bullet$\;{\em near-order} if
$\mathcal{R}$ is antisymmetric and transitive,\\
\indent\hspace{0.5cm}$\bullet$\;{\em pseudo-order} if
$\mathcal{R}$ is reflexive and antisymmetric,\\
\indent\hspace{0.5cm}$\bullet$\;{\em quasi-order} or {\em preorder}
if
$\mathcal{R}$ is reflexive and transitive,\\
\indent\hspace{0.5cm}$\bullet$\;{\em partial order} if $\mathcal{R}$
is
reflexive, antisymmetric and transitive,\\
\indent\hspace{0.5cm}$\bullet$\;{\em simple order} if $\mathcal{R}$
is
weakly complete strict order,\\
\indent\hspace{0.5cm}$\bullet$\;{\em weak order} if $\mathcal{R}$ is
complete
preorder,\\
\indent\hspace{0.5cm}$\bullet$\;{\em total order} or {\em linear
order} or {\em chain} if
$\mathcal{R}$ is complete partial order,\\
\indent\hspace{0.5cm}$\bullet$\;{\em tolerance} if $\mathcal{R}$ is
reflexive and symmetric,\\
\indent\hspace{0.5cm}$\bullet$\;{\em equivalence} if $\mathcal{R}$
is reflexive, symmetric and transitive.
\begin{rem} Clearly, universal
relation $X^2$ on a nonempty set $X$ remains a complete
equivalence relation.
\end{rem}

\noindent\textbf{Definition 7 \cite{ST}.} Let $X$ be a nonempty set
and $\mathcal{R}$ a binary relation on $X$.
\begin{enumerate}
\item [{(1)}] The inverse or transpose or dual relation of $\mathcal{R}$, denoted by $\mathcal{R}^{-1}$,
 is
defined by $\mathcal{R}^{-1}:=\{(x,y)\in X^2:(y,x)\in \mathcal{R}\}$.
\item [{(2)}] The symmetric
closure of $\mathcal{R}$, denoted by $\mathcal{R}^s$, is defined to
be the set $\mathcal{R}\cup \mathcal{R}^{-1}$
($i.e.~\mathcal{R}^s:=\mathcal{R}\cup \mathcal{R}^{-1}$). Indeed,
$\mathcal{R}^s$ is the smallest symmetric relation on $X$ containing
$\mathcal{R}$.
\end{enumerate}

\noindent{\bf Proposition 2 \cite{RT1}.} For a binary relation
$\mathcal{R}$ on a nonempty set $X$,
$$(x,y)\in\mathcal{R}^s\Longleftrightarrow [x,y]\in\mathcal{R}.$$

\noindent\textbf{Definition 8 \cite{DM}.} Let $X$ be a nonempty set,
$E\subseteq X$ and $\mathcal{R}$ a binary relation on $X$. Then, the
restriction of $\mathcal{R}$ to $E$, denoted by $\mathcal{R}|_E$, is
defined to be the set $\mathcal{R}\cap E^2$
($i.e.~\mathcal{R}|_E:=\mathcal{R}\cap E^2$). Indeed,
$\mathcal{R}|_E$ is a relation on $E$ induced by $\mathcal{R}$.\\

\noindent\textbf{Definition 9 \cite{RT1}.} Let $X$ be a nonempty set
and $\mathcal{R}$ a binary relation on $X$. A sequence
$\{x_n\}\subset X$ is called $\mathcal{R}$-preserving if
$$(x_n,x_{n+1})\in\mathcal{R}\;\;\forall~n\in \mathbb{N}_{0}.$$

\noindent\textbf{Definition 10 \cite{RT1}.} Let $X$ be a nonempty
set and $f$ a self-mapping on $X$. A binary relation $\mathcal{R}$ on $X$ is called $f$-closed if for all $x,y\in X$,
$$(x,y)\in \mathcal{R}\Rightarrow (fx,fy)\in \mathcal{R}.$$

\noindent\textbf{Definition 11.} Let $X$ be a nonempty set and $f$
and $g$ two self-mappings on $X$. A binary relation $\mathcal{R}$ on $X$ is called $(f,g)$-closed if for all $x,y\in X$,
$$(gx,gy)\in \mathcal{R}\Rightarrow (fx,fy)\in \mathcal{R}.$$

Notice that under the restriction $g=I,$ the identity mapping on
$X,$ Definition 11 reduces to Definition 10.\\

{\noindent\bf{Proposition 3.}} Let $X$ be a nonempty set,
$\mathcal{R}$ a binary relation on $X$ and $f$ and $g$ two
self-mappings on $X$. If $\mathcal{R}$ is $(f,g)$-closed, then
$\mathcal{R}^s$ is also $(f,g)$-closed.\\

In the following lines, we introduce relation-theoretic variants of the metrical notions: completeness, closedness, continuity, $g$-continuity and compatibility.\\

\noindent\textbf{Definition 12}. Let $(X,d)$ be a metric space and
$\mathcal{R}$ a binary relation on $X$. We say that $(X,d)$ is
$\mathcal{R}$-complete if every $\mathcal{R}$-preserving Cauchy
sequence in $X$ converges.

\begin{rem} Every complete metric space is
$\mathcal{R}$-complete, for any binary relation $\mathcal{R}$.
Particularly, under the universal relation the notion of
$\mathcal{R}$-completeness coincides with usual completeness.
\end{rem}

\noindent\textbf{Definition 13}. Let $(X,d)$ be a metric
space and $\mathcal{R}$ a binary relation on $X$. A subset $E$ of
$X$ is called $\mathcal{R}$-closed if every $\mathcal{R}$-preserving
convergent sequence in $E$ converges to a point of $E$.

\begin{rem}  Every closed subset of a metric
space is $\mathcal{R}$-closed, for any binary relation
$\mathcal{R}$. Particularly, under the universal relation the notion
of $\mathcal{R}$-closedness coincides with usual closedness.\end{rem}

\noindent\textbf{Proposition 4}.  An $\mathcal{R}$-complete subspace of a metric space is
$\mathcal{R}$-closed.\\
\noindent{\bf Proof.} Let $(X,d)$ be a metric space. Suppose that $Y$ is an $\mathcal{R}$-complete subspace of $X$. Take an $\mathcal{R}$-preserving sequence $\{x_n\}\subset Y$ such that $x_n\stackrel{d}{\longrightarrow} x\in X$. As each convergent sequence is Cauchy,  $\{x_n\}$ is an $\mathcal{R}$-preserving Cauchy sequence in $Y$. Hence, $\mathcal{R}$-completeness of $Y$ implies that the limit of $\{x_n\}$ must lie in $Y$, $i.e.$, $x\in Y$. Therefore, $Y$ is $\mathcal{R}$-closed.\\

\noindent\textbf{Proposition 5}. An $\mathcal{R}$-closed subspace of an
$\mathcal{R}$-complete metric space is $\mathcal{R}$-complete.\\
\noindent{\bf Proof.} Let $(X,d)$ be an
$\mathcal{R}$-complete metric space. Suppose that $Y$ is $\mathcal{R}$-closed subspace of $X$. Let $\{x_n\}$ be an $\mathcal{R}$-preserving Cauchy sequence in $Y$. As $X$ is $\mathcal{R}$-complete, $\exists~x\in X$ such that $x_n\stackrel{d}{\longrightarrow} x$ and so $\{x_n\}$ is an $\mathcal{R}$-preserving sequence converging to $x$. Hence, $\mathcal{R}$-closeness of $Y$ implies that $x\in Y$. Therefore, $Y$ is $\mathcal{R}$-complete.\\

\noindent\textbf{Definition 14.} Let $(X,d)$ be a metric space,
$\mathcal{R}$ a binary relation on $X$ and $x\in X$. A mapping
$f:X\rightarrow X$ is called $\mathcal{R}$-continuous at $x$ if for
any $\mathcal{R}$-preserving sequence $\{x_n\}$ such that
$x_n\stackrel{d}{\longrightarrow} x$, we have
$f(x_n)\stackrel{d}{\longrightarrow} f(x)$. $f$ is called
$\mathcal{R}$-continuous if it is $\mathcal{R}$-continuous at each
point of $X$.

\begin{rem} Every continuous mapping
is $\mathcal{R}$-continuous, for any binary relation $\mathcal{R}$.
Particularly, under the universal relation the notion of
$\mathcal{R}$-continuity coincides with usual continuity.\end{rem}

\noindent\textbf{Definition 15.} Let $(X,d)$ be a metric space,
$\mathcal{R}$ a binary relation on $X$, $g$ a self-mapping on $X$
and $x\in X$. A mapping $f:X\rightarrow X$ is called
$(g,\mathcal{R})$-continuous at $x$ if for any sequence $\{x_n\}$
such that $\{gx_n\}$ is $\mathcal{R}$-preserving and
$g(x_n)\stackrel{d}{\longrightarrow} g(x)$, we have
$f(x_n)\stackrel{d}{\longrightarrow} f(x)$. Moreover, $f$ is called
$(g,\mathcal{R})$-continuous if it is $(g,\mathcal{R})$-continuous at each
point of $X$.\\

Notice that under the restriction $g=I,$ the identity mapping on
$X,$ Definition 15 reduces to Definition 14.
\begin{rem} Every $g$-continuous mapping
is $(g,\mathcal{R})$-continuous, for any binary relation
$\mathcal{R}$. Particularly, under the universal relation the notion
of $(g,\mathcal{R})$-continuity coincides with usual
$g$-continuity.\end{rem}

\noindent\textbf{Definition 16.} Let $(X,d)$ be a metric space,
$\mathcal{R}$ a binary relation on $X$ and $f$ and $g$ two
self-mappings on $X$. We say that $f$ and $g$ are
$\mathcal{R}$-compatible if for any sequence $\{x_n\}\subset X$ such
that $\{fx_n\}$ and $\{gx_n\}$ are $\mathcal{R}$-preserving and
$\lim\limits_{n\to \infty}g(x_n)=\lim\limits_{n\to \infty}f(x_n)$,
we have
$$\lim\limits_{n\to \infty}d(gfx_n,fgx_n)=0.$$

\begin{rem} In a metric space $(X,d)$ endowed with a
binary relation $\mathcal{R}$, commutativity $\Rightarrow$ weak
commutativity $\Rightarrow$ compatibility $\Rightarrow$
$\mathcal{R}$-compatibility $\Rightarrow$ weak compatibility.
Particularly, under the universal relation the notion of
$\mathcal{R}$-compatibility coincides with usual compatibility
.\end{rem}

The following notion is a generalization of $d$-self-closeness of a
partial order relation $(\preceq)$ defined by Turinici
\cite{T-LCP}.\\

\noindent\textbf{Definition 17 \cite{RT1}.} Let $(X,d)$ be a metric
space. A binary relation $\mathcal{R}$ on $X$ is called
$d$-self-closed if for any $\mathcal{R}$-preserving sequence
$\{x_n\}$ such that $x_n\stackrel{d}{\longrightarrow} x$, there
exists a subsequence $\{x_{n_k}\}{\rm \;of\;} \{x_n\} \;{\rm
with}\;\;[x_{n_k},x]\in\mathcal{R}~~~\forall~k\in \mathbb{N}_{0}.$\\

\noindent\textbf{Definition 18.} Let $(X,d)$ be a metric space and
$g$ a self-mapping on $X$. A binary relation $\mathcal{R}$
on $X$ is called $(g,d)$-self-closed if for any
$\mathcal{R}$-preserving sequence $\{x_n\}$ such that
$x_n\stackrel{d}{\longrightarrow} x$, there exists a subsequence
$\{x_{n_k}\}$ of $\{x_n\}$ with
$[gx_{n_k},gx]\in\mathcal{R}~~~\forall~k\in \mathbb{N}_{0}.$\\

Notice that under the restriction $g=I,$ the identity mapping on $X,$ Definition 18 reduces to Definition 17.\\

\noindent\textbf{Definition 19 \cite{BR1}.} Let $X$ be a nonempty
set and $\mathcal{R}$ a binary relation on $X$. A subset $E$ of $X$
is called $\mathcal{R}$-directed if for each pair $x,y\in E$, there
exists $z\in X$ such that $(x,z)\in\mathcal{R}$ and
$(y,z)\in\mathcal{R}$.\\

\noindent\textbf{Definition 20 \cite{DM}.} Let $X$ be a nonempty set
and $\mathcal{R}$ a binary relation on $X$. For $x,y\in X$, a path
of length $k$ (where $k$ is a natural number) in $\mathcal{R}$ from
$x$ to $y$ is a finite sequence
$\{z_0,z_1,z_2,...,z_{k}\}\subset X$ satisfying the following conditions:\\
\indent\hspace{0.5mm} (i) $z_0=x~{\rm and}~z_k=y$,\\
\indent\hspace{0.5mm} (ii) $(z_i,z_{i+1})\in\mathcal{R}$ for each $i~(0\leq i\leq k-1)$.\\
Notice that a path of length $k$ involves $k+1$ elements of $X$,
although they are not necessarily distinct.\\

\noindent\textbf {Definition 21.} Let $X$ be a nonempty
set and $\mathcal{R}$ a binary relation on $X$. A subset $E$ of $X$
is called $\mathcal{R}$-connected if for each pair $x,y\in E$, there
exists a path in $\mathcal{R}$ from $x$ to $y$.\\

\noindent Given a binary relation $\mathcal{R}$ and two
self-mappings $f$ and $g$ defined on a nonempty set $X$, we use the
following notations:
\begin{enumerate}
\item[{(i)}] {\rm C}$(f,g):=\{x\in X:gx=fx\},~i.e.,$ the set of all coincidence points of $f$ and $g$,
\item[{(ii)}] $\overline{{\rm C}}(f,g):=\{\overline{x}\in X:\overline{x}=gx=fx,\;x\in X\},~i.e.,$ the set of all points of coincidence  of $f$ and $g$,
\item[{(iii)}] $X(f,\mathcal{R}):=\{x\in X:(x,fx)\in \mathcal{R}\}$,
\item[{(iv)}] $X(f,g,\mathcal{R}):=\{x\in X:(gx,fx)\in \mathcal{R}\}$.
\end{enumerate}

The main result of Alam and Imdad \cite{RT1} which is indeed the relation-theoretic version of Banach contraction principle runs as follows:\\

\noindent{\bf Theorem 1 \cite{RT1}.} Let $(X,d)$ be a complete metric space,
$\mathcal{R}$ a binary relation on $X$ and $f$ a self-mapping on
$X$. Suppose that the following
conditions hold:\\
\indent\hspace{0.5mm}(i) $\mathcal{R}$ is $f$-closed,\\
\indent\hspace{0.5mm}(ii) either $f$ is continuous or $\mathcal{R}$ is $d$-self-closed,\\
\indent\hspace{0.5mm}(iii) $X(f,\mathcal{R})$ is nonempty,\\
\indent\hspace{0.5mm}(iv) there exists $\alpha\in [0,1)$ such that \\
\indent\hspace{2.5cm}$d(fx,fy)\leq\alpha d(x,y)\;\;\forall~ x,y\in X$ with $(x,y)\in \mathcal{R}$.\\
Then $f$ has a fixed point. Moreover, if\\
\indent\hspace{0.5mm}(v) $X$ is $\mathcal{R}^s$-connected,\\
then $f$ has a unique fixed point.\\

Finally, we record the following known results, which are needed in the proof of our main results. \\

\noindent{\bf Lemma 1 \cite{L2}.} Let $X$ be a nonempty set and $g$
a self-mapping on $X$. Then there exists a subset $E\subseteq X$
such that $g(E)=g(X)$ and $g:E\rightarrow X$ is one-one.\\

\noindent{\bf Lemma 2 \cite{PGF13}.} Let $X$ be a nonempty set and
$f$ and $g$ two self-mappings on $X$. If $f$ and $g$ are weakly
compatible, then every point of coincidence of $f$ and $g$ is also a
coincidence point of $f$ and $g$.\\

\section{Main Results}
\label{SC:Main Results}
Now, we are equipped to prove our main result on the existence of coincidence points which runs as follows:\\

\noindent{\bf Theorem 2.} Let $(X,d)$ be a metric
space, $\mathcal{R}$ a binary relation on $X$ and $Y$ an
$\mathcal{R}$-complete subspace of $X$. Let $f$ and $g$ be two
self-mappings on $X$. Suppose that the following conditions hold:\\
\indent\hspace{5mm}$(a)$ $f(X)\subseteq g(X)\cap Y$,\\
\indent\hspace{5mm}$(b)$ $\mathcal{R}$ is $(f,g)$-closed,\\
\indent\hspace{5mm}$(c)$ $X(f,g,\mathcal{R})$ is nonempty,\\
\indent\hspace{5mm}$(d)$ there exists $\alpha\in [0,1)$ such that \\
\indent\hspace{2.5cm}$d(fx,fy)\leq\alpha d(gx,gy)\;\;\forall~ x,y\in X$ with $(gx,gy)\in \mathcal{R}$,\\
\indent\hspace{5mm}$(e)$ $(e1)$ $f$ and $g$ are $\mathcal{R}$-compatible,\\
\indent\hspace{1.2cm}$(e2)$ $g$ is $\mathcal{R}$-continuous,\\
\indent\hspace{1.2cm}$(e3)$ either $f$ is $\mathcal{R}$-continuous or $\mathcal{R}$ is $(g,d)$-self-closed,\\
\indent\hspace{2mm}or alternately\\
\indent\hspace{5mm}$(e^\prime)$ $(e^{\prime}1)$ $Y \subseteq g(X)$,\\
\indent\hspace{1.2cm}$(e^{\prime}2)$ either $f$ is $(g,\mathcal{R})$-continuous or $f$ and $g$ are continuous or $\mathcal{R}|_Y$ is $d$-self-\\
\indent\hspace{2cm} closed.\\
Then $f$ and $g$ have a coincidence point.\\

\noindent{Proof.} Assumption $(a)$ is equivalent to saying
that $f(X)\subseteq g(X)$ and $f(X)\subseteq Y$. In view of assumption $(c)$,
let $x_0$ be an arbitrary element of $X(f,g,\mathcal{R})$, then
$(gx_0,fx_0)\in \mathcal{R}$. If $g(x_0)=f(x_0)$, then $x_0$ is a
coincidence point of $f$ and $g$ and hence we are through.
Otherwise, if $g(x_0)\neq f(x_0)$, then from $f(X)\subseteq g(X)$,
we can choose $x_1\in X$ such that $g(x_1)=f(x_0)$. Again from
$f(X)\subseteq g(X)$, we can choose $x_2\in X$ such that
$g(x_2)=f(x_1)$. Continuing this process, we construct a sequence
$\{x_n\}\subset X$ (of joint iterates) such that
$$g(x_{n+1})=f(x_n)\;\;\forall ~n \in \mathbb{N}_{0}.\eqno (1)$$
Now, we claim that $\{gx_n\}$ is $\mathcal{R}$-preserving
sequence, $i.e.,$
$$(gx_n,gx_{n+1})\in \mathcal{R}\;\;\forall~ n \in \mathbb{N}_{0}.\eqno (2)$$
We prove this fact by mathematical induction. On using equation (1)
(with $n=0$) and the fact that $x_0\in X(f,g,\mathcal{R})$, we
have
$$(gx_{0},gx_1)\in \mathcal{R},$$
which shows that (2) holds for $n=0.$ Suppose that (2) holds for  $n=r>0$,
$i.e.,$
$$(gx_r,gx_{r+1})\in \mathcal{R}.$$
As $\mathcal{R}$ is $(f,g)$-closed, we have
$$(fx_r,fx_{r+1})\in \mathcal{R},$$
which, on using (1), yields that
$$(gx_{r+1},gx_{r+2})\in \mathcal{R},$$
$i.e.$, (2) holds for $n=r+1$. Hence, by induction, (2) holds for all $n \in \mathbb{N}_{0}$.\\
In view of (1) and (2), the sequence $\{fx_n\}$ is also an
$\mathcal{R}$-preserving, $i.e.,$
$$(fx_n,fx_{n+1})\in \mathcal{R}\;\;\forall~ n \in \mathbb{N}_0.\eqno (3)$$
On using (1), (2) and assumption $(d)$, we obtain
$$d(gx_{n},gx_{n+1})=d(fx_{n-1},fx_{n})\leq \alpha d(gx_{n-1},gx_{n})\;\;\;\forall~n\in \mathbb{N}.$$
By induction, we have
$$d(gx_{n},gx_{n+1})\leq \alpha d(gx_{n-1},gx_{n})\leq \alpha^2 d(gx_{n-2},gx_{n-1})\leq\cdots\leq \alpha^n d(gx_{0},gx_{1})\;\;\forall~n\in\mathbb{N}$$
so that
$$d(gx_{n},gx_{n+1})\leq \alpha^n d(gx_{0},gx_{1})\;\;\;\forall~n\in \mathbb{N}.\eqno(4)$$
For $n<m$, using (4), we obtain
\begin{eqnarray*}
d(gx_{n},gx_{m})&\leq&
d(gx_{n},gx_{n+1})+d(gx_{n+1},gx_{n+2})+\cdots+d(gx_{m-1},gx_{m})\\
&\leq& (\alpha^n+\alpha^{n+1}+\cdots+\alpha^{m-1})d(gx_{0},gx_{1})\\
&=& \frac{\alpha^n-\alpha^m}{1-\alpha}d(gx_{0},gx_{1})\\
&\leq& \frac{\alpha^n}{1-\alpha}d(gx_{0},gx_{1})\\
&\rightarrow& 0\;{\rm as}\;m,n\rightarrow \infty.
\end{eqnarray*}
Therefore $\{gx_n\}$ is a Cauchy sequence.\\

Owing to (1), $\{gx_n\}\subset f(X)\subseteq Y$ so that $\{gx_n\}$
is $\mathcal{R}$-preserving Cauchy sequence in $Y$. As $Y$ is
$\mathcal{R}$-complete, there exists $z\in Y$ such that
$$\lim\limits_{n\to\infty} g(x_n)=z.\eqno(5) $$
On using (1) and (5), we obtain
$$\lim\limits_{n\to\infty} f(x_n)=z.\eqno(6) $$

Now, we use assumptions $(e)$ and $(e^\prime)$ to accomplish the
proof. Assume that $(e)$ holds. Using (2), (5) and assumption
$(e2)$ ($i.e.$ $\mathcal{R}$-continuity of $g$), we have
$$\lim\limits_{n\to\infty} g(gx_n)=g(\lim\limits_{n\to\infty}gx_n)=g(z).\eqno(7) $$
Now, using (3), (6) and assumption $(e2)$ ($i.e.$
$\mathcal{R}$-continuity of $g$), we have
$$\lim\limits_{n\to\infty} g(fx_n)=g(\lim\limits_{n\to\infty}fx_n)=g(z).\eqno(8) $$
As $\{fx_n\}$ and $\{gx_n\}$ are $\mathcal{R}$-preserving (due to
(2) and (3)) and $\lim\limits_{n\to\infty}
f(x_n)=\lim\limits_{n\to\infty} g(x_n)=z$ (due to (5) and
(6)), on using assumption $(e1)$($i.e.$
$\mathcal{R}$-compatibility of $f$ and $g$), we obtain
$$\lim\limits_{n\to\infty}d(gfx_n,fgx_n)=0.\eqno(9)$$

Now, we show that $z$ is a coincidence point of $f$ and $g$. To
accomplish this, we use assumption $(e)$. Firstly, suppose
that $f$ is $\mathcal{R}$-continuous. On using (2), (5) and
$\mathcal{R}$-continuity of $f$, we obtain
$$\lim\limits_{n\to\infty} f(gx_n)=f(\lim\limits_{n\to\infty} gx_n)=f(z).\eqno(10) $$
On using (8), (9), (10) and continuity of $d$, we obtain
\begin{eqnarray*}
d(gz,fz)&=&d(\lim\limits_{n\to\infty}gfx_{n},\lim\limits_{n\to\infty}fgx_{n})\\
&=&\lim\limits_{n\to\infty}d(gfx_n,fgx_n)\\
&=&0
\end{eqnarray*}
so that $$ g(z)=f(z).$$ Hence we are through. Alternately, suppose
that $\mathcal{R}$ is $(g,d)$-self-closed. As $\{gx_n\}$ is $\mathcal{R}$-preserving (due
to (2)) and $g(x_n)\stackrel{d}{\longrightarrow} z$ (due to
(5)), by using $(g,d)$-self-closedness of $\mathcal{R}$, there
exists a subsequence $\{gx_{n_k}\}$ of $\{gx_n\}$ such that
$$[ggx_{n_k},gz]\in \mathcal{R}\;\;\forall~ k\in \mathbb{N}_{0}.\eqno(11)$$
Since $g(x_{n_k})\stackrel{d}{\longrightarrow} z,$ so equations (5)-(9)
hold for also $\{x_{n_k}\}$ instead of $\{x_n\}$. On using (11),
assumption $(d)$ and Proposition 1, we obtain
$$d(fgx_{n_k},fz)\leq\alpha d(ggx_{n_k},gz)\;\;\forall~ k\in \mathbb{N}_0.\eqno(12)$$
On using triangular inequality, (7), (8), (9) and (12), we get
\begin{eqnarray*}
 \nonumber d(gz,fz)&\leq& d(gz,gfx_{n_k})+d(gfx_{n_k},fgx_{n_k})+d(fgx_{n_k},fz)
\\&\leq&d(gz,gfx_{n_k})+d(gfx_{n_k},fgx_{n_k})+\alpha d(ggx_{n_k},gz)\\
&\rightarrow& 0\;{\rm as}\; k\rightarrow \infty
\end{eqnarray*}
so that $$ g(z)=f(z).$$

Thus $z\in X$ is a coincidence point of $f$ and $g$ and hence we are through.\\

Now, assume that $(e^\prime)$ holds. Owing to assumption
$(e^\prime1)$ ($i.e.$, $Y\subseteq g(X)$), we can find some $u\in X$
such that $z=g(u).$ Hence, (5) and (6) respectively reduce to
$$\lim\limits_{n\to\infty} g(x_n)=g(u).\eqno(13)$$
$$\lim\limits_{n\to\infty} f(x_n)=g(u).\eqno(14)$$

Now, we show that $u$ is a coincidence point of $f$ and $g$. To
accomplish this, we use assumption $(e^{\prime}2)$. Firstly, suppose
that $f$ is $(g,\mathcal{R})$-continuous, then using (2) and
(13), we get
$$\lim\limits_{n\to\infty} f(x_n)=f(u).\eqno(15)$$
On using (14) and (15), we get
$$g(u)=f(u).$$
Hence, we are done. Secondly, suppose that $f$ and $g$ are
continuous. Owing to Lemma 1, there exists a subset $E\subseteq X$
such that $g(E)=g(X)$ and $g:E \rightarrow X$ is one-one. Now, define $T: g(E) \rightarrow g(X)$ by
$$T(ga)=f(a)\;\;\forall\; g(a)\in g(E)\; {\rm where}\; a\in E.\eqno(16)$$
As $g:E \rightarrow X$ is one-one and $f(X)\subseteq g(X)$, $T$ is
well defined. Again since $f$ and $g$ are continuous, it follows
that $T$ is continuous. Using the fact $g(X)=g(E)$, assumptions
$(a)$ and $(e^\prime1)$ reduce to respectively $f(X)\subseteq
g(E)\cap Y$ and $Y\subseteq g(E)$, which follows that, without loss
of generality, we are able to construct $\{x_n\}_{n=1}^\infty\subset
E$ satisfying (1) and to choose $u\in E$. On using (13), (14),
(16) and continuity of $T$, we get
$$f(u)=T(gu)=T(\lim\limits_{n\to\infty} gx_n)=\lim\limits_{n\to\infty} T(gx_n)=\lim\limits_{n\to\infty} f(x_n)=g(u).$$
Thus $u\in X$ is a coincidence point of $f$ and $g$ and hence we are
through. Finally, suppose that $\mathcal{R}|_Y$ is $d$-self-closed. As $\{gx_n\}$ is $\mathcal{R}|_Y$-preserving
(due to (2)) and $g(x_n)\stackrel{d}{\longrightarrow} g(u)\in Y$
(due to (13)), using $d$-self-closedness of $\mathcal{R}|_Y$, there
exists a subsequence $\{gx_{n_k}\}$ of $\{gx_n\}$ such that
$$[gx_{n_k},gu]\in \mathcal{R}|_Y\;\;\forall~ k\in \mathbb{N}_{0}.\eqno(17)$$
On using (13), (17), assumption $(d)$ and Proposition 1, we obtain
\begin{eqnarray*}
 \nonumber d(fx_{n_k},fu)&\leq&\alpha d(gx_{n_k},gu)\\
&\rightarrow& 0~{\rm as~} k\rightarrow \infty
\end{eqnarray*}
so that
$$\lim\limits_{k\to\infty} f(x_{n_k})=f(u).\eqno(18)$$
Using (14) and (18), we get $$g(u)=f(u).$$\\
Thus, we are done. This completes the proof.\\

Now, as a consequence, we particularize Theorem 2 by assuming the $\mathcal{R}$-completeness of whole space $X$.\\

\noindent{\bf Corollary 1.} Let $X$ be a nonempty
set equipped with a binary relation $\mathcal{R}$ and a metric $d$
such that $(X,d)$ is an $\mathcal{R}$-complete metric space. Let $f$
and $g$
be two self-mappings on $X$. Suppose that the following conditions hold:\\
\indent\hspace{5mm}$(a)$ $f(X)\subseteq g(X)$,\\
\indent\hspace{5mm}$(b)$ $\mathcal{R}$ is $(f,g)$-closed,\\
\indent\hspace{5mm}$(c)$ $X(f,g,\mathcal{R})$ is nonempty,\\
\indent\hspace{5mm}$(d)$ there exists $\alpha\in [0,1)$ such that \\
\indent\hspace{2.5cm}$d(fx,fy)\leq\alpha d(gx,gy)\;\;\forall~ x,y\in X$ with $(gx,gy)\in \mathcal{R}$,\\
\indent\hspace{5mm}$(e)$ $(e1)$ $f$ and $g$ are $\mathcal{R}$-compatible,\\
\indent\hspace{1.2cm}$(e2)$ $g$ is $\mathcal{R}$-continuous,\\
\indent\hspace{1.2cm}$(e3)$ either $f$ is $\mathcal{R}$-continuous or $\mathcal{R}$ is $(g,d)$-self-closed,\\
\indent\hspace{2mm}or alternately\\
\indent\hspace{5mm}$(e^\prime)$ $(e^{\prime}1)$ there exists an $\mathcal{R}$-closed subspace $Y$ of $X$ such that $f(X)\subseteq Y\subseteq g(X)$,\\
\indent\hspace{1.2cm}$(e^{\prime}2)$ either $f$ is $(g,\mathcal{R})$-continuous or $f$ and $g$ are continuous or $\mathcal{R}|_{Y}$ is $d$-self-\\
\indent\hspace{2cm} closed.\\
Then $f$ and $g$ have a coincidence point.\\

\noindent{\bf Proof}. The result corresponding to part
$(e)$ follows easily on setting $Y=X$ in Theorem 2, while the same (result)
in the presence of part $(e^{\prime})$ follows using Proposition 5.
\begin{rem}If $g$ is onto in Corollary 1, then we can drop assumption $(a)$ as in
this case it trivially holds. Also, we can remove assumption  $(e^{\prime}1)$ as it trivially holds for $Y=g(X)=X$ using Proposition 4. Whenever, $f$ is onto, owing to assumption $(a)$, $g$ must be onto and hence again same conclusion is immediate.\end{rem}
On using Remarks 2-6, we obtain the more natural version of Theorem 2 in the form of the following consequence.\\

\noindent{\bf Corollary 2}. Theorem 2 (also Corollary 1) remains true if the usual
metrical terms namely: completeness, closedness,
compatibility (or commutativity/weak commutativity), continuity and $g$-continuity are used
instead of their respective $\mathcal{R}$-analogous.\\

Now, we present certain results enunciating the uniqueness of a point of
coincidence, coincidence point and
common fixed point corresponding to Theorem 2.\\

\noindent{\bf Theorem 3.} In addition to the hypotheses of Theorem 2, suppose that the following condition holds:\\
\indent\hspace{0.5mm} $(u_1)$: $f(X)$ is
$\mathcal{R}|_{g(X)}^s$-connected.\\
Then $f$ and $g$ have a unique point of coincidence.\\
{\noindent{Proof.}} In view of Theorem 2, $\overline{{\rm
C}}(f,g)\neq\emptyset$. Take
$\overline{x},\overline{y}\in\overline{{\rm C}}(f,g)$, then
$\exists~x,y\in X$ such that
$$\overline{x}=g(x)=f(x)\;{\rm and}\; \overline{y}=g(y)=f(y).\eqno(19)$$
Now, we show that $\overline{x}=\overline{y}.$ As $f(x),f(y)\in
f(X)\subseteq g(X)$, by assumption $(u_1)$, there exists a path (say
$\{gz_0,gz_1,gz_2,...,gz_{k}\}$) of some finite length $k$ in
$\mathcal{R}|_{g(X)}^s$ from $f(x)$ to $f(y)$ (where
$z_0,z_1,z_2,...,z_{k}\in X$). Owing to (19), without loss of
generality, we may choose $z_0=x$ and $z_k=y$. Thus, we have
$$[gz_i, gz_{i+1}]\in \mathcal{R}|_{g(X)} \;{\rm for~each}\;i\;(0\leq i\leq k-1).\eqno(20)$$
Define the constant sequences $z_n^0=x$ and $z_n^k=y$, then using
(19), we have $g(z^0_{n+1})=f(z^0_n)=\overline{x}~{\rm and~}
g(z^k_{n+1})=f(z^k_n)=\overline{y}\;\;\forall~ n\in \mathbb{N}_{0}$.
Put $z_0^1=z_1,z_0^2=z_2,..., z_0^{k-1}=z_{k-1}$. Since
$f(X)\subseteq g(X)$, on the lines similar to that of Theorem 2, we
can define sequences $\{z_n^1\},\{z_n^2\},...,\{z_n^{k-1}\}$ in $X$
such that $g(z^1_{n+1})=f(z^1_n),g(z^2_{n+1})=f(z^2_n),...,
g(z^{k-1}_{n+1})=f(z^{k-1}_n)\;\forall~ n\in \mathbb{N}_{0}$. Hence,
we have
$$g(z^i_{n+1})=f(z^i_n)\;\;\forall~ n\in \mathbb{N}_{0}\;{\rm and~for~each}\;i\;(0\leq i\leq k).\eqno(21)$$
Now, we claim that
$$[gz_n^i,gz_n^{i+1}]\in \mathcal{R}\;\;\forall~ n\in \mathbb{N}_{0}\;{\rm and~for~each}\;i\;(0\leq i\leq k-1).\eqno(22)$$
We prove this fact by the method of mathematical induction. It
follows from (20) that (22) holds for $n=0.$ Suppose that (22) holds
for  $n=r>0$, $i.e.,$
$$[gz_r^i, gz_r^{i+1}]\in \mathcal{R}\;\;{\rm for~each}\;i\;(0\leq i\leq k-1).$$
As $\mathcal{R}$ is $(f,g)$-closed, using Proposition 3, we obtain
$$[fz_r^i, fz_r^{i+1}]\in \mathcal{R}\;\;{\rm for~each}\;i\;(0\leq i\leq k-1),$$
which on using (22), gives rise
$$ [gz_{r+1}^i, gz_{r+1}^{i+1}]\in \mathcal{R}\;\;{\rm for~each}\;i\;(0\leq i\leq k-1).$$
It follows that (22) holds for $n=r+1$. Thus, by induction, (22)
holds for all $n \in \mathbb{N}_0$. Now for all $n \in \mathbb{N}_0$
and for each $i\;(0\leq i\leq k-1)$, define
$t_n^i=:d(gz_n^i,gz_n^{i+1})$. Then, we claim that
$$\lim\limits_{n\to\infty}t_n^i=0.\eqno(23)$$
On using (21), (22), assumption $(d)$ and Proposition 1, for each
$i\;(0\leq i\leq k-1)$ and for all $n \in \mathbb{N}_0$, we obtain
\begin{eqnarray*}
t_{n+1}^i&=& d(gz_{n+1}^i,gz_{n+1}^{i+1})\\
&=& d(fz_{n}^i,fz_{n}^{i+1})\\
&\leq& \alpha d(gz_{n}^i,z_{n}^{i+1})\\
&=&\alpha t_{n}^i.
\end{eqnarray*}
By induction, we have
$$t_{n+1}^i\leq \alpha t_{n}^i\leq \alpha^2 t_{n-1}^i\leq... \leq\alpha^{n+1} t_{0}^i$$
so that $$t_{n+1}^i\leq \alpha^{n+1} t_{0}^i,$$ yielding thereby
$$\lim\limits_{n\to\infty}t_{n}^i=0\;{\rm for~each}\;i\;(0\leq i\leq
k-1).$$ Thus, (23) is proved for each $i\;(0\leq i\leq k-1)$. On using
triangular inequality and (23), we obtain
$$d(\overline{x},\overline{y})\leq t_n^0+t_n^1+\cdots+t_n^{k-1}
\to 0\;\; as \;\; n\to\infty$$
$\Longrightarrow \indent\hspace{4cm}\overline{x}=\overline{y}.$\\

\noindent{\bf Corollary 3.} Theorem 3 remains
true if we replace
the condition $(u_1)$ by one of the following conditions:\\
\indent\hspace{0.5mm} $(u_1^\prime)$ $\mathcal{R}|_{f(X)}$ is complete,\\
\indent\hspace{0.5mm} $(u_1^{\prime\prime})$ $f(X)$ is
$\mathcal{R}|_{g(X)}^s$-directed.\\
{\noindent{Proof.}} If $(u_1^\prime)$ holds, then for each $u,v\in f(X)$,
$[u,v]\in\mathcal{R}|_{f(X)}\subseteq\mathcal{R}|_{g(X)}$ (owing to assumption $f(X)\subseteq g(X)$), which
amounts to say that $\{u,v\}$ is a path of length 1 in
$\mathcal{R}|_{g(X)}^s$ from $u$ to $v$. Hence $f(X)$ is
$\mathcal{R}|_{g(X)}^s$-connected consequently
Theorem 3 gives rise the conclusion.\\
Otherwise, if $(u_1^{\prime\prime})$ holds then for each $u,v\in f(X)$,
$\exists~w\in g(X)$ such that $[u,w]\in\mathcal{R}|_{g(X)}$ and
$[v,w]\in\mathcal{R}|_{g(X)}$ (owing to assumption $f(X)\subseteq g(X)$), which
amounts to say that $\{u,w,v\}$ is a path of length 2
in $\mathcal{R}|_{g(X)}^s$ from $u$ to $v$. Hence
$f(X)$ is
$\mathcal{R}|_{g(X)}^s$-connected and again by
Theorem 3 conclusion is immediate.\\

\noindent{\bf Theorem 4.} In addition to the hypotheses of Theorem 3, suppose that the following condition holds:\\
\indent\hspace{0.5mm} $(u_2)$: one of $f$ and $g$ is one-one.\\
Then $f$ and $g$ have a unique coincidence point.\\
{\noindent{Proof.}} In view of Theorem 2, ${\rm
C}(f,g)\neq\emptyset$. Take $x,y\in {\rm C}(f,g)$, then in view of
Theorem 3, we have
 $$g(x)=f(x)=f(y)=g(y).$$
 As $f$ or $g$ is one-one, we have
 $$x=y.$$

\noindent{\bf Theorem 5.} In addition to the
hypotheses embodied in condition $(e^\prime)$ of Theorem 3, suppose that the following condition holds:\\
\indent\hspace{0.5mm} $(e^\prime3)$: $f$ and $g$ are weakly compatible.\\
Then $f$ and $g$ have a unique common fixed point.\\
{\noindent{Proof.}} Owing to Remark 6 as well as
assumption $(e^\prime3)$, the mappings $f$ and $g$ are weakly compatible. Take $x\in {\rm C}(f,g)$ and denote
$g(x)=f(x)=\overline{x}$. Then in view of Lemma 2,
$\overline{x}\in {\rm C}(f,g)$. It follows from Theorem 3 with
$y=\overline{x}$ that $g(x)=g(\overline{x}),$ $i.e.$,
$\overline{x}=g(\overline{x})$, which yields that
$$\overline{x}=g(\overline{x})=f(\overline{x}).$$
Hence, $\overline{x}$ is a common fixed point of $f$ and $g$. To
prove uniqueness, assume that $x^*$ is another common fixed point of
$f$ and $g$. Then again from Theorem 3, we have
$$x^*=g(x^*)=g(\overline{x})=\overline{x}.$$
Hence we are through.\\

On setting $g = I$, the identity mapping on $X$, in Theorems 2-5, we get respectively the following corresponding fixed point
result.\\

\noindent{\bf Corollary 4.} Let $(X,d)$ be a metric
space, $\mathcal{R}$ a binary relation on $X$ and $f$ a self-mapping on
$X$. Let $Y$ be an
$\mathcal{R}$-complete subspace of $X$ such that $f(X)\subseteq Y$. Suppose that the following
conditions hold:\\
\indent\hspace{0.5mm}(i) $\mathcal{R}$ is $f$-closed,\\
\indent\hspace{0.5mm}(ii) either $f$ is $\mathcal{R}$-continuous or $\mathcal{R}|_{Y}$ is $d$-self-closed,\\
\indent\hspace{0.5mm}(iii) $X(f,\mathcal{R})$ is nonempty,\\
\indent\hspace{0.5mm}(iv) there exists $\alpha\in [0,1)$ such that \\
\indent\hspace{2.5cm}$d(fx,fy)\leq\alpha d(x,y)\;\;\forall~ x,y\in X$ with $(x,y)\in \mathcal{R}$.\\
Then $f$ has a fixed point. Moreover, if\\
\indent\hspace{0.5mm}(v) $f(X)$ is
$\mathcal{R}^s$-connected,\\
then $f$ has a unique fixed point.\\

Notice that Corollary 4 is an improvement of Theorem 1 in the following respects:\\
\noindent\vspace{.15cm} $\bullet$ Usual notions of completeness and continuity are not necessary. Alternately, they can be replaced by their respective
$\mathcal{R}$-analogues.\\
\vspace{.15cm} $\bullet$ $\mathcal{R}$-completeness of whole space $X$ and $d$-self-closedness of whole relation $\mathcal{R}$
are not necessary as they can be respectively replaced by $\mathcal{R}$-completeness
of any subspace $Y$ and $d$-self-closedness of $\mathcal{R}|_{Y}$, where$f(X)\subseteq Y\subseteq X$.\\
\vspace{.15cm} $\bullet$ For uniqueness part, $\mathcal{R}^s$-connectedness of whole space $X$ is not required
but it suffices to take the same merely of the subset $f(X)$.\\

\noindent{\bf Corollary 5.} Corollary 4 remains true if we replace
assumption (v) by one of the
following conditions:\\
\indent\hspace{0.5mm} (v)$^\prime$ $\mathcal{R}|_{f(X)}$ is complete,\\
\indent\hspace{0.5mm} (v)$^{\prime\prime}$ $f(X)$ is $\mathcal{R}^s$-directed.\\

\section{Some Consequences}
\label{SC:Some Consequences} In this section, we derive several
results of the existing literature as consequences of our newly proved
results presented in the earlier sections.\\

\subsection{Coincidence theorems in abstract metric spaces}. Under the universal relation ($i.e.$ $\mathcal{R}=X^2$), Theorems 2-5 reduce to the following coincidence point theorems:\\

\noindent{\bf Corollary 6}.  Let $(X,d)$ be a metric
space and $Y$ a complete subspace of $X$. Let $f$ and $g$ be two
self-mappings on $X$. Suppose that the following conditions hold:\\
\indent\hspace{5mm}$(a)$ $f(X)\subseteq g(X)\cap Y$,\\
\indent\hspace{5mm}$(b)$ there exists $\alpha\in [0,1)$ such that \\
\indent\hspace{2.5cm}$d(fx,fy)\leq\alpha d(gx,gy)\;\;\forall~ x,y\in X$ with $(gx,gy)\in \mathcal{R}$,\\
\indent\hspace{5mm}$(e)$ $(e1)$ $f$ and $g$ are compatible,\\
\indent\hspace{1.2cm}$(e2)$ $g$ is continuous,\\
\indent\hspace{2mm}or alternately\\
\indent\hspace{5mm}$(e^\prime)$ $Y \subseteq g(X)$.\\
Then $f$ and $g$ have a unique point of coincidence.\\

\noindent{\bf Corollary 7.} In addition to the hypotheses of Corollary 6, if
one of $f$ and $g$ is one-one, then $f$ and $g$ have a unique coincidence point.\\

\noindent{\bf Corollary 8.} In addition to the hypothesis $(e^\prime)$ of Corollary 6, if
$f$ and $g$ are weakly compatible, then $f$ and $g$ have a unique common fixed point.\\

Notice that Corollaries 6, 7 and 8 improve the well-known coincidence
theorems of Goebel \cite{CP0} and Jungck \cite{CP4}.
\subsection{Coincidence theorems under $(f,g)$-closed sets} Samet and Vetro \cite{NFI}
introduced the notion of $F$-invariant sets and utilized the same to prove
some coupled fixed point results for generalized
linear contractions on metric spaces without any partial order. Recently,
Kutbi $et\;al.$ \cite{FIC8} weakened the notion of $F$-invariant
sets by introducing the notion of $F$-closed sets. Most recently,
Karapinar $et\;al.$ \cite{FIC9} proved some unidimensional versions
of some earlier coupled fixed point results involving $F$-closed sets.
To describe such results, we need to recall the following
notions:\\

\noindent\textbf{Definition 21 \cite{FIC9}.} Let $f,g:X\rightarrow
X$ be two mappings and {\it M}$\subseteq X^2$ a subset. We
say that{\it M} is:
\begin{enumerate}
\item [{(i)}] $(f,g)$-closed if $(fx,fy)\in{\it
M}$ for all $x,y\in X$ implies that $(gx,gy)\in {\it M}$,
\item [{(ii)}] $(f,g)$-compatible if $f(x)=f(y)$ for all $x,y\in X$ implies that
$g(x)=g(y)$.\\
\end{enumerate}

\noindent\textbf{Definition 22 \cite{FIC9}.} We say that a subset
{\it M} of $X^2$ is transitive if $(x,y),(y,z)\in {\it M}$ implies
that $(x,z)\in {\it M}$.\\

\noindent\textbf{Definition 23 \cite{FIC9}.} Let $(X, d)$ be a
metric space and {\it M}$\subseteq X^2$ a subset. We say that
$(X,d,{\it M})$ is regular if for all sequence $\{x_n\}\subseteq X$
such that $x_n\stackrel{d}{\longrightarrow} x$ and $(x_n, x_{n+1})\in {\it M}$ for all
$n$, we have $(x_n,x)\in {\it M}$ for all $n$.\\

\noindent\textbf{Definition 24 \cite{FIC9}.} Let $(X,d)$ be a metric
space, {\it M}$\subseteq X^2$ a subset and $x\in X$. A
mapping $f:X\rightarrow X$ is said to be {\it M}-continuous at $x$
if for all sequence $\{x_n\}\subseteq X$ such that $x_n\stackrel{d}{\longrightarrow}
x$ and $(x_n, x_{n+1})\in {\it M}$ for all $n$, we have
$f(x_n)\stackrel{d}{\longrightarrow} f(x)$. Moreover, $f$ is called {\it M}-continuous if it is {\it
M}-continuous at each
$x\in X$.\\

The following notion is introduced in order to improve the
commutativity condition of the pair of mappings $f$ and $g$, which is inspired by
the notion of $O$-compatibility of Luong and Thuan \cite{CP5} in ordered
metric spaces.\\

\noindent\textbf{Definition 25 \cite{FIC9}.}  Let $(X,d)$ be a
metric space and {\it M}$\subseteq X^2$. Two mappings
$f,g:X\rightarrow X$ are said to be {\it M}-compatible if
$$\lim\limits_{n\to \infty}d(gfx_n,fgx_n)=0$$
whenever $\{x_n\}$ is a sequence in $X$ such that $(gx_n,
gx_{n+1})\in {\it M}$ for all $n$ and $\lim\limits_{n\to
\infty}f(x_n)=\lim\limits_{n\to \infty}g(x_n)\in X$.\\

Notice that Karapinar $et\;al.$ \cite{FIC9} (inspired by the notion
of $O$-compatibility in \cite{CP5}) preferred to call``{\it(O,M)}-compatible" instead of ``{\it M}-compatible". Here the
symbol "$O$" has no pertinence as Luong and Thuan \cite{CP5} used the
term ``$O$-compatible" due to available partial ordering on the underlying metric space
($i.e.$ $O$ means order relation). But in above context, Karapinar
$et\;al.$ \cite{FIC9} used a nonempty subset {\it M} without partial ordering, so it is appropriate to use the term ``{\it
M}-compatible".\\

Here, it can be point out that the involved set {\it M} being a
subset of $X^2$ is indeed a binary relation on $X$. Therefore, the
concept of $(f,g)$-closed subset of $X^2$ can be interpreted as
$(f,g)$-closed binary relation on $X$. Obviously, Definitions 16
and 17 are weaker than Definitions 25 and 23 respectively. Taking
$\mathcal{R}=${\it M} in Corollary 1, we get an improved version of the following result of Karapinar
$et\;al.$ \cite{FIC9}.\\

\noindent{\bf Corollary 9} (see Corollary 34 \cite{FIC9}).  Let
$(X,d)$ be a complete metric space, let $f,g:X\rightarrow X$ be two
mappings and let {\it M}$\subseteq X^2$
be a subset such that\\
\indent\hspace{0.5mm}(i) $f(X)\subseteq g(X)$,\\
\indent\hspace{0.5mm}(ii) {\it M} is $(f,g)$-compatible and $(f,g)$-closed,\\
\indent\hspace{0.5mm}(iii) there exists $x_{0}\in X$ such that $(gx_{0},fx_{0})\in {\it M}$,\\
\indent\hspace{0.5mm}(iv) there exits $\alpha\in [0,1)$ such
that\\
\indent\hspace{2.5cm}$d(fx,fy)\leq\alpha d(gx,gy)\;\;\forall~ x,y\in X$ with $(gx,gy)\in {\it M}$.\\
Also assume that, at least, one of the following conditions holds:\\
\indent\hspace{0.5mm}$(a)$ $f$ and $g$ are {\it M}-continuous and
{\it M}-compatible,\\
\indent\hspace{0.5mm}$(b)$ $f$ and $g$ are continuous and commuting,\\
\indent\hspace{0.5mm}$(c)$ $(X,d,{\it M})$ is regular and $g(X)$ is closed.\\
Then $f$ and $g$ have, at least, a coincidence point.\\

Observe that {\it M}-compatibility of $(f,g)$(see assumption (ii)) is unnecessary.\\

\subsection{Coincidence theorems in ordered metric spaces via increasing mappings} Indeed the present trend was initiated by Turinici \cite{P15,P16}, Ran and Reurings
\cite{PF1} and Nieto and Rodr\'{\i}guez-L\'{o}pez \cite{PF2} which
was was later generalized by many authors (e.g.
\cite{PGF2,PGF3,PGF13}). In this subsection as well as in succeeding
subsection, $X$ denotes a nonempty set endowed with a partial order
$\preceq$. In what follows, we write $\succeq:=\preceq^{-1}$ and
$\prec\succ:=\preceq^{s}$. On the lines of O'Regan and Petru\c{s}el
\cite{PGF2}, the triple $(X,d,\preceq)$ is called ordered metric
space wherein $X$ denotes a nonempty set endowed with
a metric $d$ and a partial order $\preceq$.\\

\noindent\textbf{Definition 26 \cite{PGF3}.} Let $(X,\preceq)$ be an
ordered set and  $f$ and $g$ two self-mappings on $X$. We say that
$f$ is $g$-increasing if for any $x,y\in X$, $g(x)\preceq
g(y)\Rightarrow f(x)\preceq f(y)$.
\begin{rem} It is clear that  $f$ is $g$-increasing iff $\preceq$ is $(f,g)$-closed.\end{rem}
\noindent\textbf{Definition 27 \cite{PGF13}.} Given a mapping
$g:X\rightarrow X$, we say that an ordered metric space
$(X,d,\preceq)$ has {\it g-ICU}\;(increasing-convergence-upper
bound) property if $g$-image of every increasing sequence $\{x_n\}$
in $X$ such that $x_n\stackrel{d}{\longrightarrow} x$, is bounded
above by $g$-image of its limit (as an upper bound), $i.e.,$
$g(x_n)\preceq g(x)\;\;\forall~ n\in \mathbb{N}_{0}.$\\
Notice that under the restriction $g=I,$ the identity mapping on
$X,$ Definition 27 transforms to the notion of {\it ICU} property.
\begin{rem} It is clear that if $(X,d,\preceq)$ has {\it ICU} property (resp. {\it g-ICU} property), then $\preceq$ is $d$-self-closed (resp. $(g,d)$-self-closed).\end{rem}

On taking $\mathcal{R}=\preceq$ in Corollary 2 and using Remarks
8 and 9, we obtain the following result, which is an improved version of
Corollary 3 of Alam $et\;al.$ \cite{PGF13}.\\

\noindent{\bf Corollary 10.} Let $(X,d,\preceq)$ be an
ordered metric space and $Y$ a complete subspace
of $X$. Let
$f$ and $g$ be two self-mappings on $X$. Suppose that the following conditions hold:\\
\indent\hspace{5mm}$(a)$ $f(X)\subseteq g(X)\cap Y$,\\
\indent\hspace{5mm}$(b)$ $f$ is $g$-increasing,\\
\indent\hspace{5mm}$(c)$ there exists $x_{0}\in X$ such that $g(x_{0})\preceq f(x_{0})$,\\
\indent\hspace{5mm}$(d)$ there exists $\alpha\in [0,1)$ such that \\
\indent\hspace{2.5cm}$d(fx,fy)\leq\alpha d(gx,gy)\;\;\forall~ x,y\in X$ with $g(x)\preceq g(y)$,\\
\indent\hspace{5mm}$(e)$ $(e1)$ $f$ and $g$ are compatible,\\
\indent\hspace{1.2cm}$(e2)$ $g$ is continuous,\\
\indent\hspace{1.2cm}$(e3)$ either $f$ is continuous or $(Y,d,\preceq)$ has {\it g-ICU} property,\\
\indent\hspace{6mm}or alternately\\
\indent\hspace{5mm}$(e^\prime)$ $(e^{\prime}1)$ $Y \subseteq g(X)$,\\
\indent\hspace{1.2cm}$(e^{\prime}2)$ either $f$ is
$g$-continuous or $f$ and $g$ are continuous or
$(Y,d,\preceq)$ has\\
\indent\hspace{2cm} {\it ICU} property.\\
Then $f$ and $g$ have a coincidence point.\\

\subsection{Coincidence points in ordered metric spaces via comparable mappings} The core results involving comparable mappings are contained in Nieto and Rodr\'{\i}guez-L\'{o}pez \cite{PF3}, Turinici
\cite{T-RRF,T-NLF}, Dori\'{c} $et\;al.$ \cite{PF-C0} and Alam and
Imdad \cite{PGF15}.\\

\noindent\textbf{Definition 28 \cite{PGF15}.} Let $(X,\preceq)$ be
an ordered set and  $f$ and $g$ two self-mappings on $X$. We say
that $f$ is $g$-comparable if for any $x,y\in X,$
$$g(x)\prec\succ g(y)\Rightarrow f(x)\prec\succ f(y).$$
\begin{rem} It is clear that  $f$ is $g$-comparable iff $\prec\succ$ is $(f,g)$-closed.\end{rem}

\noindent\textbf{Definition 29 \cite{PGF14}.} Let $(X,\preceq)$ be
an ordered set and $\{x_n\}\subset X.$
\begin{enumerate}
\item[{(i)}] the sequence $\{x_n\}$ is said to be termwise bounded if there is an element $z\in X$ such that each term
of $\{x_n\}$ is comparable with $z,$ $i.e.$,
$$x_n\prec\succ z\;\;\;\;\;\;\;\forall~ n\in \mathbb{N}_0$$
so that $z$ is a c-bound of $\{x_n\}$ and
\item[{(ii)}] the sequence $\{x_n\}$ is said to termwise monotone if consecutive terms of $\{x_n\}$ are comparable,
$i.e.$,
$$x_n\prec\succ x_{n+1}\;\;\forall~ n\in \mathbb{N}_0.$$
\end{enumerate}

\begin{rem} Clearly, $\{x_n\}$ is termwise monotone iff it is $\prec\succ$-preserving.
\end{rem}

\noindent\textbf{Definition 30 \cite{PGF14}.} Given a mapping
$g:X\rightarrow X$, we say that an ordered metric space
$(X,d,\preceq)$ has {\it g-TCC}\;(termwise
monotone-convergence-c-bound) property if every termwise monotone
sequence $\{x_n\}$ in $X$ such that
$x_n\stackrel{d}{\longrightarrow} x$ has a subsequence,  whose
$g$-image is termwise bounded by $g$-image of limit (of the
sequence) as a c-bound, $i.e.,$ $g(x_{n_k})\prec\succ g(x)\;\forall~
k\in \mathbb{N}_{0}.$\\
Notice that under the restriction $g=I,$ the identity mapping on $X,$ Definition 30 transforms to the notion of {\it
TCC} property.
\begin{rem} Clearly, $(X,d,\preceq)$ has {\it TCC} property (resp. {\it g-TCC} property) iff $\prec\succ$ is $d$-self-closed (resp. $(g,d)$-self-closed).\end{rem}

On taking $\mathcal{R}=\prec\succ$ in Corollary 2 and using Remarks
10 and 12, we obtain the following result, which is an improved
version of Theorem 3.7 of Alam and Imdad \cite{PGF15}.\\

\noindent{\bf Corollary 11.} Let $(X,d,\preceq)$ be an
ordered metric space and $Y$ a complete subspace
of $X$. Let
$f$ and $g$ be two self-mappings on $X$. Suppose that the following conditions hold:\\
\indent\hspace{5mm}$(a)$ $f(X)\subseteq g(X)\cap Y$,\\
\indent\hspace{5mm}$(b)$ $f$ is $g$-comparable,\\
\indent\hspace{5mm}$(c)$ there exists $x_{0}\in X$ such that $g(x_{0})\prec\succ f(x_{0})$,\\
\indent\hspace{5mm}$(d)$ there exists $\alpha\in [0,1)$ such that \\
\indent\hspace{2.5cm}$d(fx,fy)\leq\alpha d(gx,gy)\;\;\forall~ x,y\in X$ with $g(x)\prec\succ g(y)$,\\
\indent\hspace{5mm}$(e)$ $(e1)$ $f$ and $g$ are compatible,\\
\indent\hspace{1.2cm}$(e2)$ $g$ is continuous,\\
\indent\hspace{1.2cm}$(e3)$ either $f$ is continuous or $(Y,d,\preceq)$ has {\it g-TCC} property,\\
\indent\hspace{6mm}or alternately\\
\indent\hspace{5mm}$(e^\prime)$ $(e^{\prime}1)$ $Y \subseteq g(X)$,\\
\indent\hspace{1.2cm}$(e^{\prime}2)$ either $f$ is
$g$-continuous or $f$ and $g$ are continuous or
$(Y,d,\preceq)$ has\\
\indent\hspace{2cm} {\it TCC} property.\\
Then $f$ and $g$ have a coincidence point.\\
\subsection{Coincidence theorems under symmetric
closure of a binary relation} The origin of such results can be
traced back to Samet and Turinici \cite{BR1} which is also pursued
in Berzig \cite{BR2}. In this context, $\mathcal{R}$ stands for an
arbitrary
 binary relation on a nonempty set $X$ and $\mathcal{S}:=\mathcal{R}^s$.\\

\noindent\textbf{Definition 31 \cite{BR2}.} Let $f$ and $g$ be two
self-mappings on $X$. We say that $f$ is $g$-comparative if for any
$x,y\in X,$
$$(gx,gy)\in \mathcal{S}\Rightarrow (fx,fy)\in \mathcal{S}.$$
\begin{rem} It is clear that  $f$ is $g$-comparative iff $\mathcal{S}$ is $(f,g)$-closed.\end{rem}

\noindent\textbf{Definition 32 \cite{BR1}.} We say that
$(X,d,\mathcal{S})$ is regular if the following condition holds: if
the sequence $\{x_n\}$ in $X$ and the point $x\in X$ are such that
$$(x_n,x_{n+1})\in \mathcal{S}\;{\rm for~all~} n\; {\rm and~}\lim\limits_{n\to\infty} d(x_n,x)=0,$$
then there exists a subsequence $\{x_{n_k}\}$ of $\{x_n\}$ such that
$(x_{n_k},x)\in \mathcal{S}$ for all $k$.
\begin{rem} Clearly, $(X,d,\mathcal{S})$ is regular iff $\mathcal{S}$ is $d$-self-closed.\end{rem}

Taking the symmetric closure $\mathcal{S}$ of an arbitrary relation
$\mathcal{R}$ in Corollary 2 and using Remarks 13 and 14, we obtain
an improved version of the
following result of Berzig \cite{BR2}.\\

\noindent{\bf Corollary 12} (see Corollary 4.5 \cite{BR2}). Let $(X,d)$ be a metric space, $\mathcal{R}$ a binary relation on $X$ and $f$ and $g$ two self-mappings on $X$. Suppose that the following conditions hold:\\
\indent\hspace{0.5mm}$(a)$ $f(X)\subseteq g(X)$,\\
\indent\hspace{0.5mm}$(b)$ $f$ is $g$-comparative,\\
\indent\hspace{0.5mm}$(c)$ there exists $x_0\in X$ such that $(gx_0,fx_0)\in \mathcal{S}$,\\
\indent\hspace{0.5mm}$(d)$ there exists $\alpha\in [0,1)$ such that \\
\indent\hspace{2.5cm}$d(fx,fy)\leq\alpha d(gx,gy)\;\;\forall~ x,y\in X$ with $(gx,gy)\in \mathcal{S}$,\\
\indent\hspace{0.5mm} $(e)$ $(X,d)$ is complete and $g(X)$ is closed,\\
\indent\hspace{0.5mm} $(f)$ $(X,d,\mathcal{S})$ is regular.\\
 Then $f$ and $g$ have a coincidence point.\\

\section{Examples}
\label{SC:Examples} In this section, we provide two examples establishing the utility of
Theorems 2-5.
\begin{example} Consider $X=\mathbb{R}$ equipped with usual metric and also define a binary relation $\mathcal{R}=\{(x,y)\in \mathbb{R}^2:|x|-|y|\geq 0\}$. Then $(X,d)$ is an $\mathcal{R}$-complete metric space. Consider the mappings $f,g:X\to X$ defined by
$f(x)=\frac{x^2}{3}\;{\rm and}\;g(x)=\frac{x^2}{2}\;\forall~ x\in
X.$ Clearly, $\mathcal{R}$ is $(f,g)$-closed. Now, for $x,y\in X$ with $(gx,gy)\in \mathcal{R}$, we
have
$$d(fx,fy)=\left|\frac{x^2}{3}-\frac{y^2}{3}\right|=\frac{2}{3}\left|\frac{x^2}{2}-\frac{y^2}{2}\right|=\frac{2}{3}d(gx,gy)< \frac{3}{4}d(gx,gy).$$
Thus, $f$ and $g$ satisfy assumption $(d)$ of Theorem 2 with
$\alpha=\frac{3}{4}$. By a routine calculation, one can verify all the conditions
mentioned in $(e)$ of Theorem 2. Hence all the conditions of Theorem
2 are satisfied for $Y=X$, which guarantees that $f$ and $g$ have a coincidence point in $X$.
Moreover, observe that $(u_1)$ holds and henceforth in view of Theorem 3,
$f$ and $g$ have a unique point of coincidence (namely: $\overline{x}=0$), which remains also a unique common fixed point (in view of Theorem 5).\\
Observe that the underlying binary relation $\mathcal{R}$ is a
preorder which is not antisymmetric and henceforth not a partial
order. Thus, in all, our results are genuine extension of several corresponding results proved under partial ordering.
\end{example}
\begin{example} Consider $X=\mathbb{R}$ equipped with usual metric and also define a binary relation $\mathcal{R}=\{(x,y)\in \mathbb{R}^2:x\geq 0,\;
 y\in \mathbb{Q}\}$. Consider the mappings $f,g:X\to X$ defined by $f(x)=1\;{\rm and}\;g(x)=x^2-3\;\forall~ x\in
X.$ Clearly, $\mathcal{R}$ is $(f,g)$-closed. Now, for $x,y\in X$ with $(gx,gy)\in \mathcal{R}$, we
have
$$d(fx,fy)=|1-1|=0\leq\alpha|x^2-y^2|=\alpha d(gx,gy).$$
Thus, $f$ and $g$ satisfy assumption $(d)$ of Theorem 2 for any arbitrary
$\alpha\in [0,1)$. Also, the mappings $f$ and $g$ are not $\mathcal{R}$-compatible and hence
$(e)$ does not hold. But the subspace $Y:=g(X)=[-3,\infty)$ is
$\mathcal{R}$-complete and $f$ and $g$ are continuous, $i.e.$, all
the conditions mentioned in $(e^{\prime})$ are satisfied. Hence, in view of
Theorem 2, $f$ and $g$ have a coincidence point in $X$. Further,
in this example $(u_1)$ holds and henceforth, owing to Theorem
3, $f$ and $g$ have a unique point of coincidence (namely:
$\overline{x}=1$). Notice that neither $f$ nor $g$ is one-one,
$i.e.$, $(u_2)$ does not hold and hence, we can not apply Theorem
4, which guarantees the uniqueness of coincidence point. Observe
that, in the present example, there are two coincidence points (namely: x=2 and x=-2). Also,
$f$ and $g$ are not weakly compatible , $i.e.$, $(e^\prime3)$ does not
hold and hence, we can not apply Theorem 5, which ensures the
uniqueness of common fixed point. Notice that there is no common fixed point of $f$ and $g$.\\
Observe that the underlying binary relation $\mathcal{R}$ is a transitive relation. Indeed, $\mathcal{R}$ is non-reflexive, non-irreflexive, non-symmetric
as well as non-antisymmetric and hence it is not a preorder, partial
order, near order, strict order, tolerance or equivalence and also never turns out to be a
symmetric closure of any binary relation.
\end{example}
Here, it can be point out that corresponding results
contained in Section 4 cannot be used in the context of
present example, which substantiate the
utility of our newly proved coincidence theorems over corresponding several relevant results.

\end{document}